\documentclass[mca,article,submit,moreauthors,pdftex]{Definitions/mdpi} 
\pdfoutput=1

\firstpage{1} 
\pubvolume{xx}
\issuenum{1}
\articlenumber{5}
\pubyear{2019}
\copyrightyear{2019}
\history{Received: date; Accepted: date; Published: date}
\preto{\abstractkeywords}{\nolinenumbers}

\numberwithin{equation}{section}

\usepackage{float}
\usepackage[T1]{fontenc}
\usepackage[utf8]{inputenc}
\usepackage{epstopdf}
\usepackage{cite}
\usepackage[justification=centering]{caption}
\usepackage{tabu}
\usepackage{amssymb}
\usepackage[font=small,labelfont=bf]{caption}
\usepackage{wrapfig}
\usepackage{caption}
\usepackage{subcaption}

\usepackage{amsmath}
\usepackage{dblfloatfix}
\usepackage{graphicx}
\usepackage{color}

\newcommand{\be}{\begin{equation}}
\newcommand{\ee}{\end{equation}}

\Title{Jacob's Ladder: Prime numbers in 2d}

\Author{Alberto Fraile $^{1,}$*, Roberto Mart\'inez $^{2}$ and Daniel Fern\'andez $^{3}$}

\AuthorNames{Alberto Fraile, Roberto Mart\'inez and Daniel Fern\'andez}

\address{%
$^{1}$ \quad Czech Technical University Prague. Department of Control Engineering, Advanced Materials Group Praha 2, Karlovo náměstí 13, E-s136, Czech Republic; albertofrailegarcia@gmail.com\\
$^{2}$ \quad Universidad del País Vasco - Euskal Herriko Unibertsitatea, Barrio Sarriena s/n 48940 Leioa. Bizkaia, Spain; rrmartinezz@yahoo.es\\
$^{3}$ \quad University of Iceland, Science Institute, Dunhaga 3, 107 Reykjav\'ik, Iceland; fernandez@hi.is}

\corres{Correspondence: albertofrailegarcia@gmail.com}

\abstract{Prime numbers are one of the most intriguing figures in mathematics. Despite centuries of research, many questions remain still unsolved. In recent years, computer simulations are playing a fundamental role in the study of an immense variety of problems. In this work, we present a simple representation of prime numbers in two dimensions that allows us to formulate a number of conjectures that may lead to important avenues in the field of research on prime numbers. In particular, although the zeroes in our representation grow in a somewhat erratic, hardly predictable way, the gaps between them present a remarkable and intriguing property: a clear exponential decay in the frequency of gaps vs gap size. The smaller the gaps, the more frequently they appear. Additionally, the sequence of zeroes, despite being non-consecutive numbers, contains a number of primes approximately equal to $n/\log n$ , being $n$ the number of terms in the sequence.}


\begin{document}

\section{Introduction}
\label{sec:Introduction}

Prime numbers have fascinated mathematicians since the beginnings of Mathematics (\citeauthor{ingham, littlewood,anglin, apostol}). Their distribution is intriguing and mysterious to some extent, despite being non-chaotic. After centuries of research there are still many open problems to be solved, and the exact details of the prime numbers distribution are yet to be understood (\citeauthor{hardy} (\citeyear{hardy}), \citeauthor{shanks}). Today, the interest in prime numbers has received a new impulse after their unexpected appearance in different contexts ranging from cryptology (\citeauthor{stallings}) or quantum chaos (\citeauthor{goles, toha}) to biology (\citeauthor{berry, sakhr}). Furthermore, with the advent of more and more powerful computers, different studies are being undertaken where methods traditionally used by physicists are being applied to the study of primes. For example, in (\citeauthor{wolf} (\citeyear{wolf})) the multifractality of primes was investigated, whereas some appropriately defined Lyapunov exponents for the distribution of primes were calculated numerically in (\citeauthor{gamba}). Additionally, one of the most studied problems in number theory, the Goldbach conjecture, first formulated around 1740, has been being continuously studied with the help of computers (\citeauthor{saouter, richstein, granville})\footnote{Its weaker version, called the ternary Goldbach conjecture, was finally proven to be true in 2013 (\citeauthor{helfgott}).}. More examples can be found in (\citeauthor{borwein}).

\section{Jacob's Ladder}
\label{sec:Jacob}

One can argue that prime numbers present perplexing features, somewhat of a hybrid of local unpredictability and global regular behavior. It is this interplay between randomness and regularity what motivated searches for local and global patterns that could perhaps be signatures of more fundamental mathematical properties.

Our work is concerned with the long standing question of the prime numbers distribution, or more precisely with the gaps between primes (\citeauthor{cramer, erdos}), a topic that has attracted much attention recently (\citeauthor{granville, maynard1} (\citeyear{maynard1})) after some massive advances (\citeauthor{zhang}). 
Some problems related to prime gaps are well known, for instance, Legendre's conjecture states that there is a prime number between $n^2$ and $(n + 1)^2$ for every positive integer $n$. This conjecture is one of Landau's problems on prime numbers (\citeauthor{guy}) and up to date (2018) it is still waiting to be proved or disproved.
Another example would be the famous twin primes conjecture (\citeauthor{hardy2} (\citeyear{hardy2})). Proving this conjecture seemed to be far out of reach until just recently. In 2013, Yitang Zhang demonstrated the existence of a finite bound $B = 70,000,000$, such that there are infinitely many pairs of distinct primes which differ by no more than $B$ (\citeauthor{zhang}). After such an important breakthrough, proving the twin prime conjecture looked much more plausible. Immediately after, a cooperative team led by Terence Tao, building upon Zhang's techniques, was able to lower that bound down to 4,680. Then in the same year, Maynard slashed the value of $B$ down to 600, and finally the Polymath project further reduced it to 246 (\citeauthor{polymath1}, \citeauthor{maynard3} (\citeyear{maynard3})).

Displaying numbers in two dimensions has been a traditional approach toward primes visualization (\citeauthor{stein, orlowski}). Here we propose an original way of number arrangement yielding an appealing visual structure: an oscillating plot that increases and decreases according to the prime number distribution. We plot the integers from 1 to $n$ in 2D ($x$, $y$) starting with 1 and $y = 0$ (hence, the first point will be (1, 0)), moving to 2 the plot moves up in the y axis, so that the next point in terms of coordinates is (2, 1). In the next step, it goes up or down depending on whether the number $n$ is a prime number or not. Number 2 is prime, so it flips in such a way that it goes down next, and hence the third point is (3, 0). Now, 3 is prime so the next step goes up again, and we move up to (4, 1). Number 4 is not a prime so it continues moving up, and so on. We used our own code but different codes to produce Jacob’s ladder can be found in the On-Line Encyclopedia of Integer Sequences (\citeauthor{link}).

\begin{figure}[h]
	\centering
	\begin{minipage}{.49\textwidth}
		\centering
		\includegraphics[width=.96\linewidth]{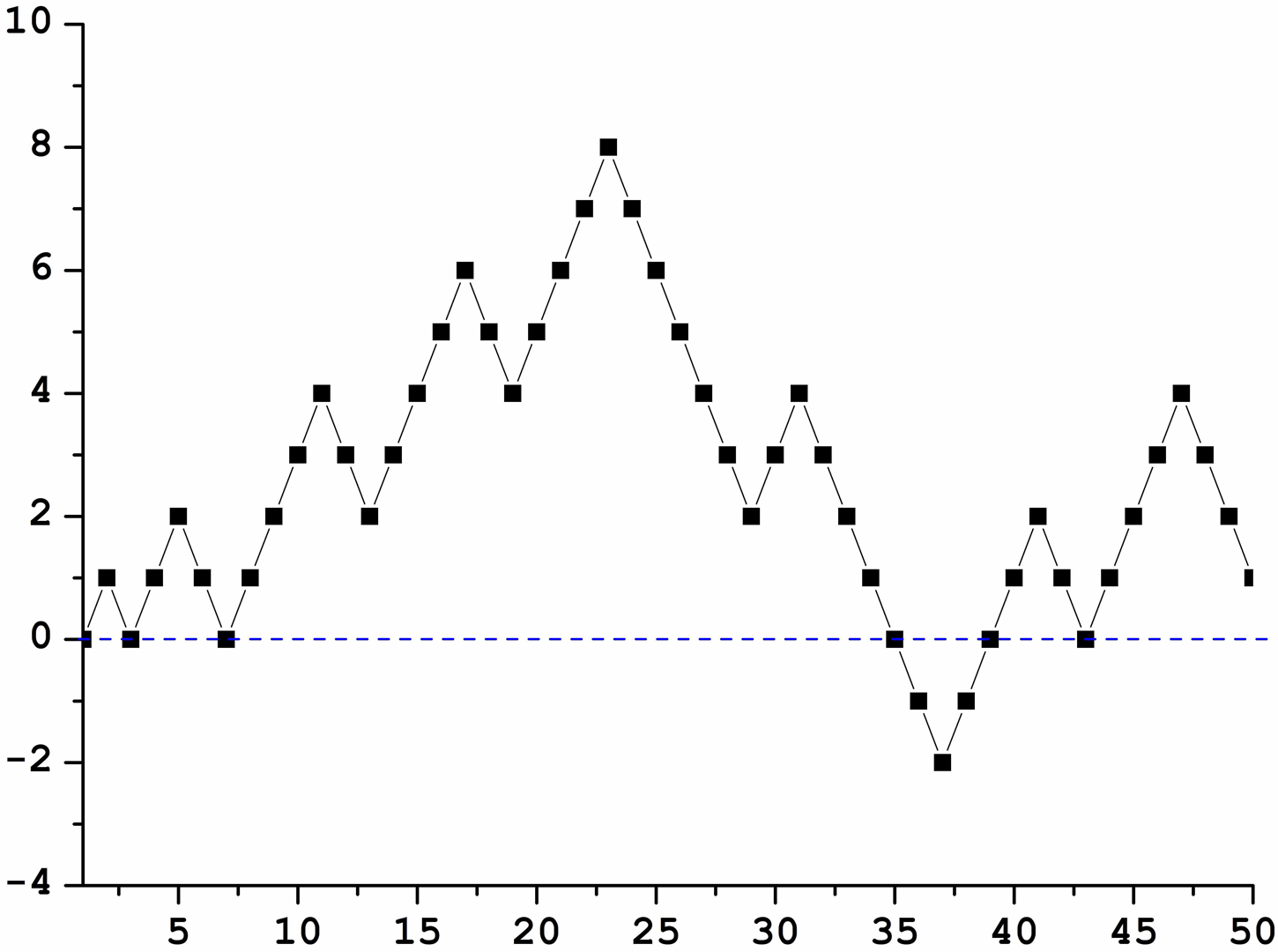}
		\caption{Illustrative plot of the first 50 points of the Jacob's Ladder sequence. The blue dashed line stands for $y = 0$.}
		\label{fig:first}
	\end{minipage}
	\begin{minipage}{.49\textwidth}
		\centering
		\includegraphics[width=.96\linewidth]{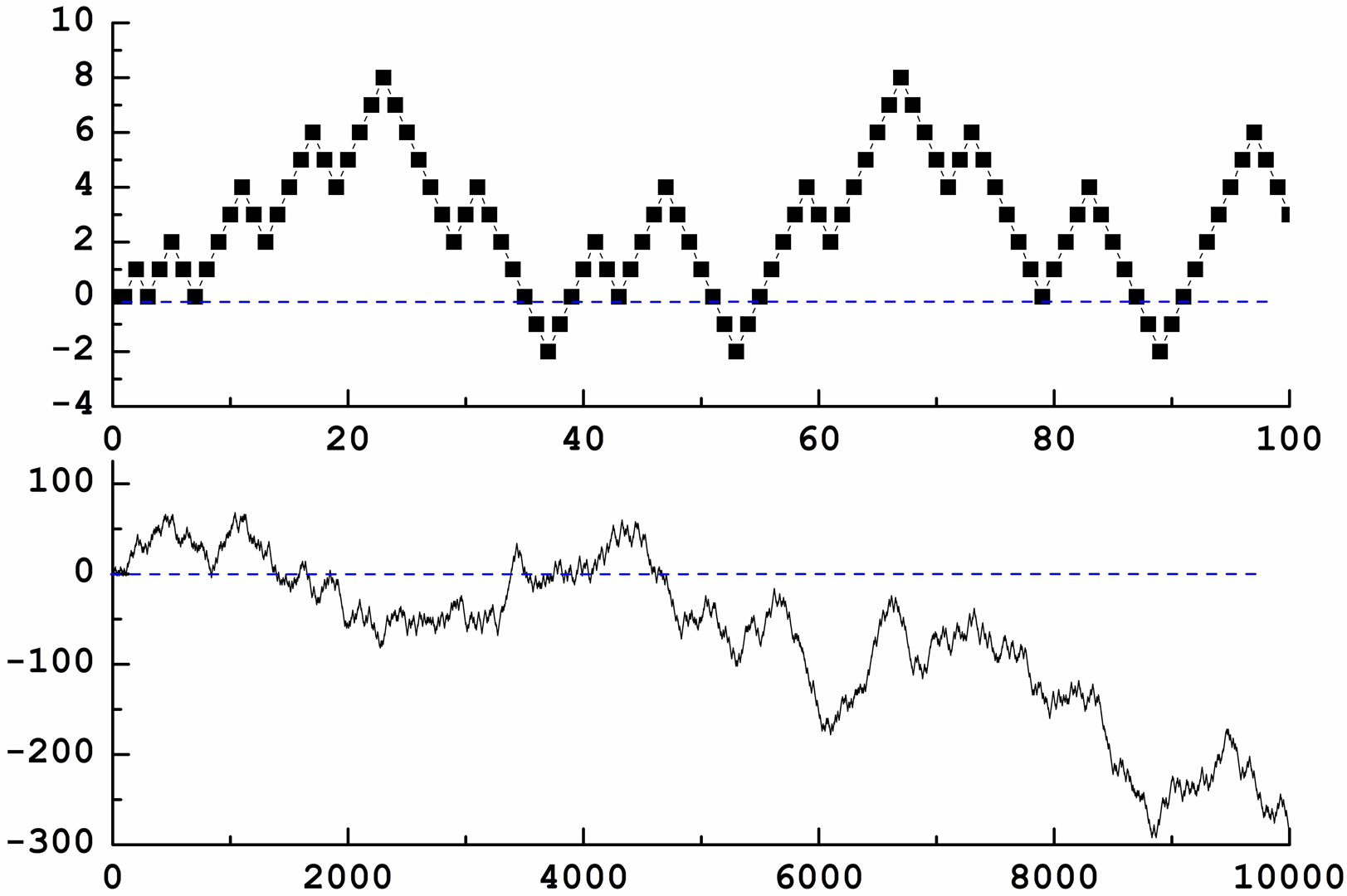}
		\caption{\emph{Top:} Jacob’s Ladder from 1 to 100. \emph{Bottom:} Jacob’s Ladder from 1 to 10,000. The blue dashed line stands for $y=0$.}
		\label{fig:second}
	\end{minipage}
\end{figure}

In Fig. \ref{fig:first} the sequence produced by the algorithm is shown up to $n=50$. The blue dashed line stands for the $y = 0$ line (or $x$ axis), which will be central to our study. We will refer to the points ($x$, 0) as zeroes from now on. Because of the resemblance of this numerical structure to a ladder, we refer to this set of points in the 2D plane as \emph{Jacob’s Ladder}, $J(n)$ for short hereafter. The points in the ladder (the $y$ values) can be written as:
\be
J(n) = \sum_{k=1}^{n} (-1)^{\pi(k)}
\ee
where $\pi(n)$ denotes the number of primes $\leq n$.
 
\section{Results}
\label{sec:Results}

The next figures show the numerical results which inspired the ideas and conjectures we will present in the next section. Fig. \ref{fig:second} (Top) shows Jacob's Ladder from 1 to 100. The blue dashed line signals the $x$ axis for clarity shake. As can be seen, up to 100 the Ladder is almost positive, being most of its points above the $x$ axis. This is misleading, though, as we can see from looking at the behavior of Jacob's Ladder from 1 to 10,000 (Fig. \ref{fig:second} (Bottom)).

\begin{wrapfigure}{L}{0.5\textwidth}
	\centering
	\includegraphics[width=.98\linewidth]{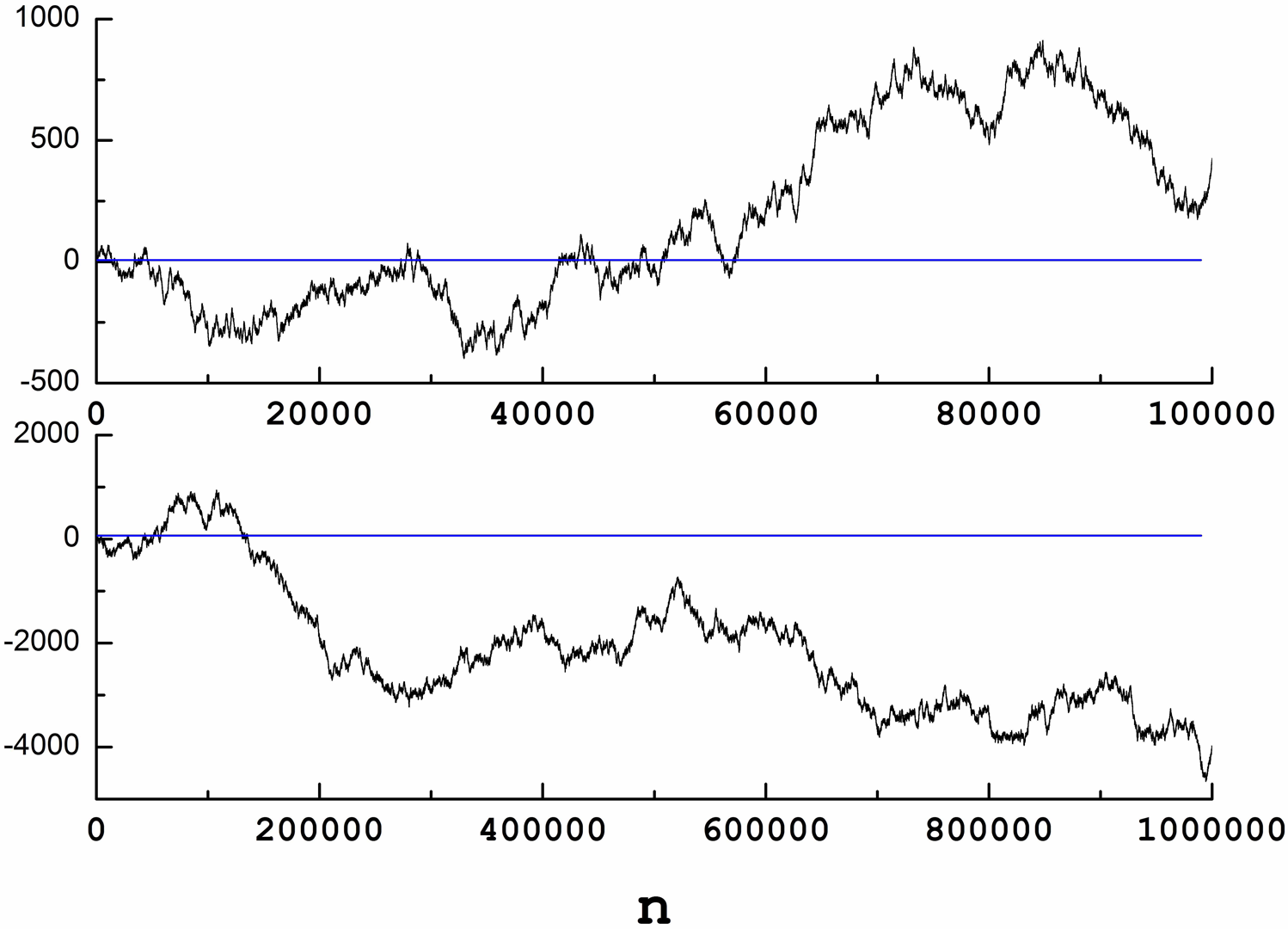}
	\caption{Jacob’s Ladder from 1 to 100,000 (top) and from 1 to 1,000,000 (bottom). The blue line stands for $y=0$.}
	\label{fig:third}
\end{wrapfigure}

\begin{figure}[!b]
	\centering
	\includegraphics[width=0.72\textwidth]{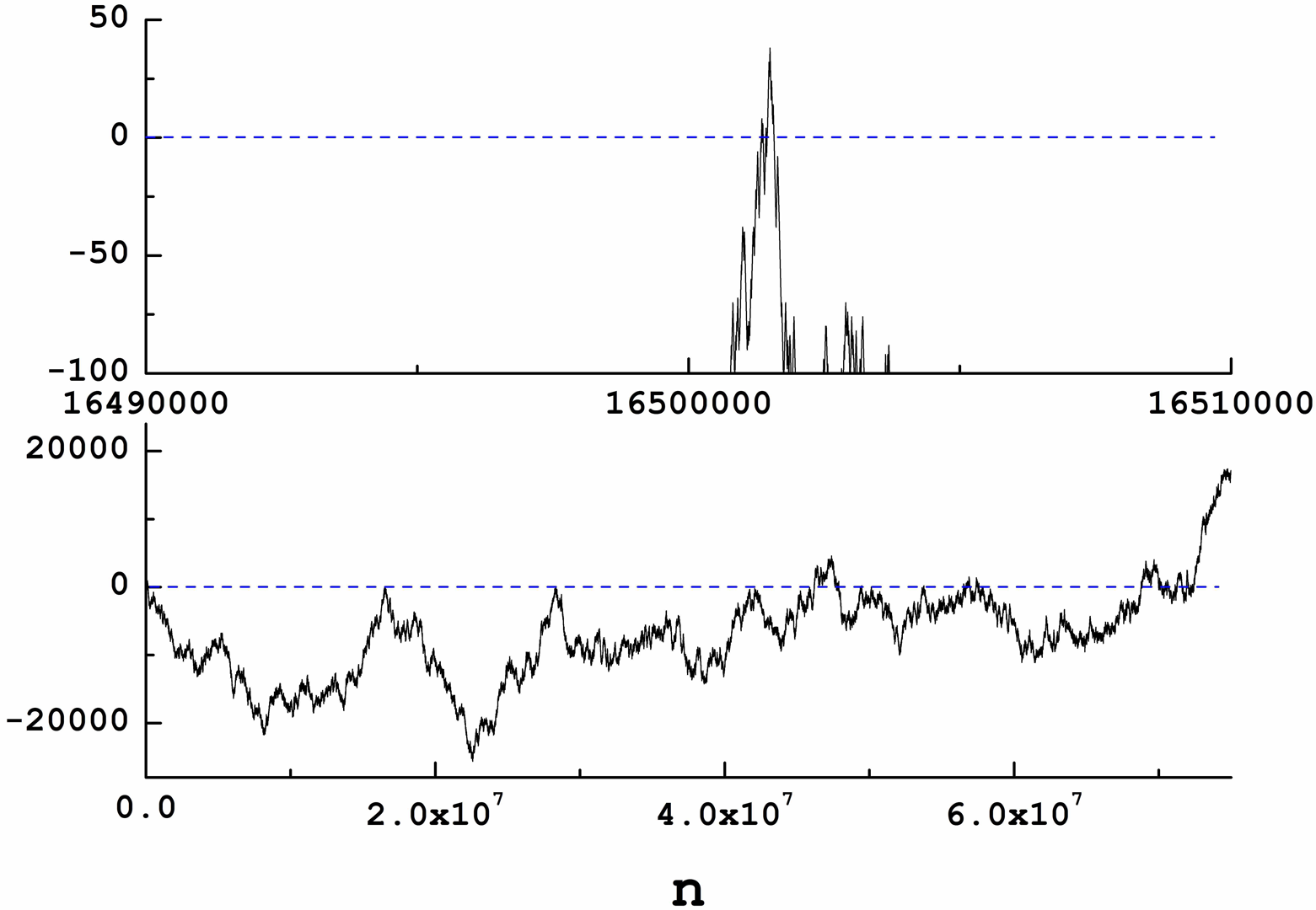} \\
	\caption{\emph{Top:} Jacob’s Ladder goes up again at $\sim 16,000,000$. \emph{Bottom:} Jacob’s Ladder shows a big positive peak and hundreds of zeroes around 44 million (and around 70 million). The blue dashed line stands for $y=0$.}
	\label{fig:fourth}
\end{figure}

Fig. \ref{fig:second} (bottom) shows Jacob’s Ladder from 1 to 10,000. As can be seen, now most of the Ladder is negative except for two regions. However, this fact changes again if we move to higher values. Fig. \ref{fig:third} (top) shows Jacob’s Ladder from 1 to 100,000 and here we see, as said before, that after $n\sim 50,000$ most of the Ladder is again positive. Fig. \ref{fig:third} (bottom) shows Jacob's Ladder from 1 to 1,000,000. Note that the Ladder presents a big region of negative values after $n\sim 150,000$. After this, no more zeroes are present up until 16,500,000 (See Fig. \ref{fig:fourth} (Top)). Afterwards, the Ladder is mostly negative again. 

Around 45 million, hundreds of zeroes are found (Fig. \ref{fig:fourth}, (bottom)). Going up to 100 million, more zeroes appear, totaling 2,415. The next crossing appears for $n=202,640,007$. Larger gaps appear in the sequence at this point. For instance, no zeroes are found between 6,572,751,127 and 9,107,331,715.

\vspace{6pt}

\noindent Much larger is the gap found at 9,169,470,319. The next zero is 222,564,950,675. And an even larger is the one from 225,888,362,295 to 676,924,242,679 (which are consecutive zeroes as well).

Jacob's ladder resembles what in probability theory is called a Wiener process, so the probability distributions could be studied as well. It is quite clear that points very far from the $x$ axis will be less abundant when compared to the number of points close to it, above or below the axis, so renormalizing the number of points in the Ladder up to a given point, it seems natural to expect some version of an inverse arcsine law (centered around the $x$ axis) in the distribution of these points. Despite the large interval explored, the data are not enough to plot a convincing graph due to the ladder's noticeable asymmetry, so we prefer to leave this point as an open question.

Building from the figures presented above, this simple representation of prime numbers in two dimensions brings about many interesting questions and allows us to arrive at a list of conjectures, which, despite their simplicity, will unfortunately be very difficult to prove or disprove. The most natural one is discussed below. Two additional ones are presented and analyzed in Appendix 1. The continuous appearance of an increasing number of zeroes, although sparse, points towards all of them holding true in the limit of large numbers.

\vspace{12pt}

{\bf Conjecture I:}

The number of cuts (zeroes) in the $x$ axis tends to infinity.
In mathematical language, being $Z(n)$ the number of zeroes in the Ladder,
\be
\lim^{}_{n\to\infty} Z(n) = \infty
\ee

\underline{Discussion}:

Proving this conjecture is beyond the scope of this article, but in the following we describe the empirical motivation behind it. This is the main idea in our study and the basis of the conjectures that follow. Fig. \ref{fig:fifth} presents the number of zeroes vs $n$, in Jacob’s Ladder, from 1 to $8\cdot 10^{11}$ (logarithmic scale on both axes). Notice that the number of zeroes increases (by construction, it cannot decrease) with $n$ in an apparently chaotic or unpredictable way. In some intervals of $n$ the number of zeroes is constant, meaning that the ladder is entirely above or below $y=0$. After those plateaus, it increases again and again. If our conjecture is true, it will increase forever as we move towards increasingly bigger values of $n$. However, we conjecture as well that the slope dictating this increase will be lower and lower as $n$ goes to infinity. The fact that the prime numbers become increasingly separated seems to indicate so: since the prime numbers are more and more separated, the ladder will present a smaller amount of zigzagging.

The sequence of cuts, or zeroes, can be denoted as $F_0(n)$ and in the same way infinitely many successions can be defined: $F_1(n)$ being the number of cuts on the $y=1$ line up to $n$, $F_{-4}(n)$ being the number of cuts in $y=-4$ up to $n$, and so on. In other words, we conjecture that $F_0(n)$ contains infinitely many elements, and, if Conjecture I is true, then it is likely that any other $F_y(n)$ will contain infinitely many terms when $n$ tends to infinity. Thus, the Ladder allows to non-trivially define infinitely many successions with infinitely many terms (and without repetitions), whose cardinal will be $\aleph_0$ and the sum of all of them is again $\aleph_0$.

It is to be noted that the ladder cannot be constrained between two $y$ values, since the gaps between primes can be arbitrarily large. A lower or upper limit could exist, but not both.

An interesting question can be formulated here: how does $Z(n)$ grow? What would be a good approximation for the number of zeroes given a value of $n$? The answer is not trivial, since the $Z(n)$ function is neither multiplicative nor additive. In Fig. \ref{fig:fifth} we present two simple functions that operate as approximate upper and lower limits of the $Z(n)$ function, namely $\sqrt{n}$ and $\sqrt[3]{n}$, represented by red lines. The number of zeroes is not strictly confined between these two curves, however the agreement is fair enough for the purposes of this paper. We also plotted the counting function $\pi(n)$ in blue, for the shake of comparison.

The idea here is not to extract accurate upper and lower bounds with these functions, but to simply offer a qualitative approximation. The reason why such a qualitative conclusion is relevant lies on the geometrical insight that can be extracted from these simple approximations. If the numbers from 1 to $n^2$ were to be plotted within a square (for instance forming a spiral, like in (\citeauthor{stein})), then having approximately $n$ zeroes (or any other selection of numbers used as tracking mechanism) would mean that the amount of zeroes we have is approximately given by the length of the side of the square. In an analogous way, if we display the numbers from 1 to $n^3$ within a cube of side $n$, then the side of the cube would give us that information. Therefore, having a number of zeroes constrained between $\sqrt{n}$ and $\sqrt[3]{n}$ seems to point towards some kind of fractality for the number of zeroes present in the sequence $F_0(n)$.

\section{Discussion}
\label{sec:Discussion}

While completing this paper, we found similar studies, like those presented in (\citeauthor{billingsley, peng, wolf2} (\citeyear{wolf2})). For instance, in (\citeauthor{wolf2} (\citeyear{wolf2})) a one-dimensional random walk (RW), where steps up and down are performed according to the occurrence of special primes (twins and cousins), was defined. If there are infinitely many twins and cousins (as suggested by the Hardy-Littlewood conjecture), then the RW defined there will continue to perform steps forever, in contrast to the RW considered in (\citeauthor{billingsley}) or (\citeauthor{peng}), where random walks were finite. In our case, we prefer not to talk about random walks since the distribution of primes, despite mysterious, is not random. On the other side, the idea of Jacob’s ladder is simple and beautiful given that, due to the infinitude of prime numbers, its intrinsic complexity and zigzagging must continue indefinitely. Now, its properties, by definition, depend both on the number of primes in a given interval and on the separation between them.
It is known that the gaps between consecutive prime numbers cluster on multiples of 6. Because of this fact 6 is sometimes called the jumping champion, and it is conjectured that it stays the champion all the way up to about $10^{35}$ (\citeauthor{wolf3} (\citeyear{wolf3}), \citeauthor{odlyzko}).
\begin{figure}[!t]
	\centering
	\includegraphics[width=0.72\textwidth]{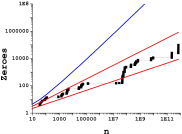} \\
	\caption{(Color online) $Z(n)$, the number of zeroes in Jacob's Ladder from 1 to $8\cdot 10^{11}$. Note the logarithmic scale on both axes. Number of zeroes = 194,531. Red lines stand for the $\sqrt{n}$ and $\sqrt[3]{n}$ functions. In blue, we plotted the counting function, $\pi(n)$, for the sake of comparison.}
	\label{fig:fifth}
\end{figure}

Beyond $10^{35}$, and until $10^{425}$, the jumping champion becomes 30 (=$2 \cdot 3 \cdot 5$), and beyond that point, the most frequent gap is 210 (=$2 \cdot 3 \cdot 5 \cdot 7$) (\citeauthor{odlyzko}). It is a natural conjecture that after some large number, the jumping champion will be 2,310 (=$2 \cdot 3 \cdot 5 \cdot 7 \cdot 11$) and so on. Further interesting results on some statistical properties between gaps have been recently found (\citeauthor{szpiro1} (\citeyear{szpiro1}), \citeauthor{szpiro2} (\citeyear{szpiro2})). However, all the aforementioned numerical observations, despite revealing intriguing properties of the primes sequence, are not easily put into our problem in order to know whether or not the Ladder will have infinitely many zeroes or not. On the other side, according to Ares et al. (\citeauthor{ares}), the apparent regularities observed in some works (\citeauthor{wolf2} (\citeyear{wolf2}), \citeauthor{wolf4} (\citeyear{wolf4}), \citeauthor{ball}) reveal no structure in the sequence of primes, and that is precisely a consequence of its randomness (\citeauthor{hull}). However, this is a controversial topic. Recent computational work points that ``after appropriate rescaling, the statistics of spacings between adjacent prime numbers follows the Poisson distribution''. See (\citeauthor{wolf5} (\citeyear{wolf5}), \citeauthor{garcia}) and references therein for more on the statistics of the gaps between consecutive prime numbers. 

The question at hand is encapsulated by whether after a given number all turns in the Ladder will have an ``order''. By order we mean that $J(n)$ will go up (or down) after the next prime, and then down (or up) and so on but always keeping a property: that an ``averaged'' curve will continue asymptotically along one particular direction of the plane without ever going back to the $x$ axis. Or in other words, is it possible that, after an unknown number $X$, the sum of all intervals between primes ``up'' minus the sum of all intervals ``down'' (or the other way around) will be positive (or negative) for all possible $Y>X$ values as $n\to\infty$?

\be
\displaystyle\sum_{X}^{Y} \,(I[\uparrow]-I[\downarrow])>0\,(<0)\;\forall Y>X
\label{eq:II}
\ee

\noindent where $I[\uparrow]$ ($I[\downarrow]$) stands for the interval between two primes [$p_n-p_{n-1}$] which make to Ladder go up (down)\footnote{\emph{Sensu stricto} this could be true for a number $n=X$ and the sum represented by Eq. \ref{eq:II} could be $> z$ (or $< z$) without crossing the $x$ axis because we start from some point up (or down) the axis.}. That looks unlikely and would be some interesting order property if it turned out to be the case (and it would prove Conjecture I to be false). If so, it would be interesting to find that number $X$ distinctly.
\begin{figure}[!t]
	\centering
	\includegraphics[width=0.72\textwidth]{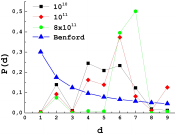} \\
	\caption{(Color online) Leading digit histogram of the zeroes sequence $F_0(n)$. The proportion of the different leading digits found in the intervals, $[1, 10^{10}]$, $[1, 10^{11}]$ and in $[1, 8\cdot 10^{11}]$ is shown in black squares, red diamonds and green circles as labeled. The expected values according to Benford's law are shown in blue triangles.}
	\label{fig:eighth}
\end{figure}

This is related to the question of whether the zeroes are randomly distributed, or follow some kind of ordered distribution. For instance, do the terms in $F_0(n)$ follow Benford's law (\citeauthor{benford})? In many different natural datasets and mathematical sequences, the leading digit\footnote{The leading digit of a number is its non-zero leftmost digit. For example, the leading digits of 2,018 and 994,455 are 2 and 9 respectively.} $d$ is found to be non-uniformly distributed as opposed to naive expectations. Instead, it presents a biased probability given by $P(d) = \log_{10} (1+1/d)$ where $d=1, 2, \dots 9$. While this empirical law was indeed first discovered by the Canadian–American astronomer, applied mathematician and autodidactic polymath Simon Newcomb (\citeauthor{newcomb}) in 1881, it is popularly known as Benford's law, and also as the law of anomalous numbers (\citeauthor{benford}).
Up to date, many mathematical sequences such as $(n^n)_{n\in \mathbb{N}}$ and $(n!)_{n\in \mathbb{N}}$ (\citeauthor{klafter}), binomial arrays ${n}\choose{k}$ (\citeauthor{draconis}), geometric sequences or sequences generated by recurrence relations (\citeauthor{raimi, takloo}), to cite a few, have been proved to conform to Benford.
\begin{figure}[!t]
	\centering
	\includegraphics[width=0.72\textwidth]{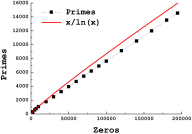} \\
	\caption{(Color online) Number of primes as a function of the number of terms (zeroes) in the sequence $F_0(n)$. The red line stands for the plausible similitude with the counting function $\pi(n) \sim n/\log n$.}
	\label{fig:ninth}
\end{figure}

Fig. \ref{fig:eighth} shows the proportion of the different leading digits found in the different intervals $[1,n]$, up to 194,531 zeroes, and the expected values according to Benford's law. Note that intervals $[1, n]$ have been chosen such that $n=10^D$, with $D$ a natural integer, in order to ensure an unbiased sample where all possible first digits are \emph{a priori} equiprobable. With that in mind, however, the results for the interval $[1,\, 8\cdot 10^{11}]$ are also presented because the number of zeroes are many times more than those found up to $10^{11}$. Obviously, the curve is likely to change as more and more zeroes are counted (it changes when fewer zeroes are counted). A straightforward calculation demonstrates that, in the best of cases, one should count up to five (or six) times the current number of zeroes before having a distribution in fair agreement with Benford's law. Note that in order to obtain such large amount of zeroes, according to Fig. \ref{fig:fifth}, it is likely that one should explore the Ladder in a 25 (36) or 125 (216) times larger interval, i.e. up to $10^{14}$ ($1.7\cdot 10^{14}$).

On the other hand, it is natural to assume that the distribution presented by $F_0(n)$ will be the same for all values of $y$. Thus, while being aware that the numbers are still very small from a statistical point of view, the analysis of the terms in more crowded sequences can provide further evidence for the random distribution (or better said, non-Benford-like distribution) of $F_y(n)$ for all $y$.

Another point to be noticed is that all of the zeroes can be written as $(4n+3)$. This is not a surprise, in fact it can be easily proved that it must be this way: after 3, which is the first zero excluding 1, all zeroes must be separated by a number twice the sum of the number of gaps between primes (every pair of primes is separated by an even number), therefore the gaps between zeroes will be always a multiple of 4.

For one more test about some possible unexpected non-random properties of the zeroes in the Ladder we can check the number of zeroes ending in a given odd digit (since all the even ones will not be in $F_0(n)$). Doing so we observe they change depending on the explored interval, so there is no clear evidence of any trend \footnote{The terms ending in 1, 3, 5, 7 and 9 amount to 38,835, 38,905, 38,799, 38,898 and 39,093 respectively in the first 194,530 zeroes, thus showing a uniform distribution ($\sim 20\%$ each). A linear fit ($y=a+bx$) to the data (in probabilities, given as $\%$) provides the following result: $a = 19.93459$ (0.04315) and $b = 0.01308$ (0.00751). The slope is small, and just twice the error of the linear fit, so it is quite unsafe to assume a trend on this.}.

But going back to the zeroes as a whole, we can observe the two most important results found in our study. First, it is natural to wonder: how many of these zeroes are primes? Is this number somewhere predictable? Since in an interval $[1, n]$ we find $\sim n/ \log n$ primes, should we expect, in a set of $X$ consecutive zeroes (but not consecutive numbers!), an amount of primes approximately equal to $X/\log X$? Interestingly, it seems to be the case if we count a few thousands of zeroes when looking at Fig. \ref{fig:ninth}, but it is not really the case from looking at the differences in percentage (see Table \ref{tab:first}). The prime zeroes found are nevertheless in remarkable agreement with this simple assumption. As we see in Fig. \ref{fig:ninth}, for $n > 10^9$ the prime numbers found in $F_0(n)$ depart steadily from $n/\log n$, being the latter an overestimate of the actual number of primes.

As mentioned before, all zeroes can be written as $(4n+3)$. But $4n+3$ is the same as $4(n+1)-1$. It is a well known result that all primes (except 2) can be written as $4n+1$ or $4n-1$ 
- the first ones allowing to be expressed as the sum of two perfect squares $p=x^2 +y^2$, while the second ones can never be expressed in such a way. Therefore in our results, those zeroes which are prime always belong to the second subgroup.
Furthermore, the arithmetic progressions $4n\pm 1$, (and $6n \pm 1$) are of particular interest because the primes in the progression $4n-1$ (and $6n-1$) are known to be "denser" than in the progression $4n+1$ (and $6n+1$); a result conjectured in 1853 by Tschebyschef (\citeauthor{new1,new2}) and proved by Phragmen (\citeauthor{new3}) and Landau (\citeauthor{new4}). However, the problem is, as is usually the case in number theory, far from being simple. A famous result by Littlewood (\citeauthor{littlewood}) showed that there are infinitely many integers $x$ such that there are more prime numbers in the progression $4n+1$ than in the $4n-1$, up to that number $x$ \footnote{Using formal mathematical notation, it is said that $\pi_{4,3}(x) < \pi_{4,1}(x)$ for infinitely many $x$, where $\pi_{b,c}(x)$ denotes the number of positive primes $\leq x$ which are equal to c (mod b).}.

It is easy to see that every sequence $F_b(x)$ (or zeroes in the lines $y = b$) will contain primes (that is, reversals of the ladder) if and only if $b$ is even. In particular if $b = 4a$, with $a$ an integer number or zero, then those primes will belong to the arithmetic progression $4n-1$. Besides, those sequences where $b = 4a \pm 2$ will contain primes belonging to the arithmetic progression $4n+1$ (with the exception of 2, there are no primes in the $F_b(x)$ sequences when $b$ is odd). Thus, it is clear that the problem of how to demonstrate the existence of an infinite number of zeroes in Jacob’s Ladder is related with previous important results regarding primes in arithmetic progressions. However, it is not enough to know how many primes are found in a particular progression: the gaps between primes are equally important, and that is a more mysterious problem.

\begin{figure}[!t]
	\centering
	\includegraphics[width=0.72\textwidth]{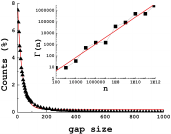} \\
	\caption{Number of times that a gap appears depending on the the gap size, only up to 1000 for clarity shake. The red line presents the exponential decay fit. The inset shows the averaged gap $\Gamma$ vs the interval size.}
	\label{fig:tenth}
\end{figure}

Now we will see what can be said about gaps: how are the zeroes separated? As we have seen in the Figures \ref{fig:first} to \ref{fig:fourth}, the Ladder goes up and down in an erratic way, resulting in a rather capricious distribution of the zeroes. Nevertheless, counting the number of times that a given gap between zeroes appears within a given number of zeroes, the result is clearly non-arbitrary. Fig. \ref{fig:tenth} shows the result after analyzing the gaps in $F_0(n)$ up to $8\cdot 10^{11}$. It can be seen that a clear exponential decay is found. Note as well that all gaps are multiples of 4, except 2, which only appears once (The demonstration of the fact that all gaps are multiples of 4 is trivial).

The inset shows the averaged gap, $\Gamma$, vs the interval size. $\Gamma$ is calculated as the sum of gaps divided by the number of gaps. Another possible definition, $\Xi$, would be the interval size divided by the number of zeroes. The plateaus are due to intervals (from $10^6$ to $10^7$ for instance) where no more zeroes are found, so the number of gaps is constant. Using the $\Xi$ definition, the averaged gap would be 10 times larger.

Fig. \ref{fig:tenth}, showing the distribution of gaps, is a remarkable result and it may lead to a number of interesting observations. Note that almost 25\% of the gaps are 4, 8, 12, or 16, and more than 53\% of the gaps are smaller than 60.
Gaps $\leqslant 100$ represent more than 63\% of the gaps found, however, as said before, very large gaps are found as well
as large as 451,035,880,384 (about half the interval explored).

In the intervals explored up to now, 4 is always the most frequent gap, referred to as $\gamma_1$ hereafter, followed by $\gamma_2$= 8, $\gamma_3$= 12 (i.e. the second most frequent gap is 8, the third most frequent is 12, etc.). Immediately, a legitimate and important question would be to inquire about the behavior of the distribution of gaps for large values of $n$. Will $\gamma_1$ be 4 for any interval, no matter how large? Or is it possible that from a very large interval onwards, it shifts to 8, and later on to 12, and so on?
Furthermore, it could be reasonable as well to expect an equiprobable distribution of gaps in the limit $n \to\infty$.

However, it has been shown that 4 is the most probable gap size (Fig. \ref{fig:tenth}). If the exponential decay of $\gamma_n$ vs $n$ hypothesis were finally confirmed $\forall n$, such result would represent an attribute of the Ladder of fundamental importance, meaning that every time that the ladder crosses the $x$ axis, it is close to a prime number. This observation can be misleading, though, because at the same time the average gap clearly seems to increase with $n$ (Inset in Fig. \ref{fig:tenth}).
So these results do not mean that all of the zeroes are primes. In fact, as showed before (See Fig. \ref{fig:ninth} and Table \ref{tab:first}), out of $n$ zeroes, about $n/\log n$ are actually primes.

\begin{table}
	\begin{tabu} to \textwidth { | X[l] | X[l] | X[l] | X[l] | X[l] | X[l] | }
		\hline \vspace{0.03cm}
		Interval & \vspace{0.03cm} Zeroes & \vspace{0.03cm} Primes & \vspace{0.03cm} $n/\log n$ & \vspace{0.03cm} Diff ($\%$) & \vspace{0.03cm} Average gap $\Gamma$ \vspace{0.06cm} \\
		\hline
		$[1, 10^2]$ & 10 & 5 & 4.342 & 13.16 & 9.2 \\
		\hline
		$[1, 10^3]$ & 16 & 6 & 5.770 & 3.82 & 9.25 \\
		\hline
		$[1, 10^4]$ & 59 & 21 & 14.469 & 31.09 & 36.20 \\
		\hline
		$[1, 10^5]$ & 139 & 36 & 28.169 & 21.75 & 526.57 \\
		\hline
		$[1, 10^6]$ & 151 & 37 & 30.096 & 18.65 & 1503.97 \\
		\hline
		$[1, 10^7]$ & 151 & 37 & 30.096 & 18.65 & 1503.97 \\
		\hline
		$[1, 10^8]$ & 2,415 & 313 & 310.034 & 0.947 & 40170.11 \\
		\hline
		$[1, 10^9]$ & 7,730 & 846 & 863.41 & -2.058 & 887722.55 \\
		\hline
		$[1, 10^{10}]$ & 11,631 & 1,161 & 1,242.438 & -7.014 & 523588.07 \\
		\hline
		$[1, 10^{11}]$ & 11,631 & 1,161 & 1,242.438 & -7.014 & 523588.07 \\
		\hline
		$[1, 8\cdot 10^{11}]$ & 194,530 & 14,556 & 15,973 & -9.734 & 2750072.04 \\
		\hline
	\end{tabu}
	\caption{Summary of the results insofar. The first column is the interval of $J(n)$, the second is the number of zeroes found in it, the third is the number of primes in the $F_0(n)$ sequence, the fourth is $n/\log(n)$ being $n$ the number of zeroes, the fifth shows the difference in percentage, and the sixth the average gap between zeroes.}
	\label{tab:first} 
\end{table}

\begin{figure}[!t]
	\centering
	\includegraphics[width=0.72\textwidth]{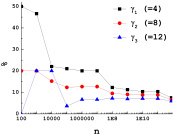} \\
	\caption{Percentage of the total gaps of value 4, 8 and 12 (colors black, red and blue respectively) vs number of zeroes, in log scale. As can be seen the percentage, decreases for all these three cases, with the order remaining as $\gamma_1$=4, $\gamma_2$=8, $\gamma_3$=12.}
	\label{fig:eleventh}
\end{figure}

The percentage of gaps of length 4 decreases with $n$ faster than exponentially, in the same way as those of length 8 or 12, and likely for all gap values (in small intervals some variations may occur, but we are interested in properties appearing in the limit of large intervals), as shown in Fig. \ref{fig:eleventh}. This result can be easily explained since the larger the interval is, the more gaps enter in the list of gaps, so it is natural that the percentage would decrease.

At the same time, since the average gap increases, it could be natural to expect an increase in the percentage of larger gaps. But that is not the case -- only for small intervals do we observe an increase with $n$ of the percentage of large gaps, because for small intervals large gaps are rare, yet it seems that the percentage of large gaps is always below that of 4, 8 etc.

A final important remark: we started this work without knowing any previous studies on this topic. Only while writing down the paper that we checked if the sequence of zeroes appeared in the On-Line Encyclopedia of Integer Sequences (\citeauthor{link}), and then we found it was first studied in 2001 by Jason Earls and subsequently by Hans Havermann within a different context \footnote{See full sequence in \href{http://chesswanks.com/num/a064940.txt}{http://chesswanks.com/num/a064940.txt}.} and Don Reble \footnote{See full sequence in \href{https://oeis.org/A064940/a064940.txt}{https://oeis.org/A064940/a064940.txt}.}.

\section{Conclusions}
\label{sec:Conclusions}

In this paper, an original idea has been proposed with the aim of arranging integers in 2D according to the occurrence of primality. This arrangement of numbers resembles a `ladder' and displays peaks and valleys at the positions of the primes. Numerical studies have been undertaken in order to extract qualitative information from such a representation. As a result of those studies, it is possible to promote four observations in the category of conjectures. One may easily state that the ladder, from its own construction, will cross the horizontal axis an indefinite number of times; namely, the ladder will have infinitely many zeroes.

The number of zeroes grows in a hardly predictable way, somewhere between $\sqrt{n}$ and $\sqrt[3]{n}$, but the number of primes in $F_0(n)$ is very well fitted by $n/\log n$ for small numbers (up to 10,000 zeroes). However, for larger numbers the fit is just fair. Interestingly, this sequence of zeroes does not seem to follow Benford's Law, at least up to the range of $n = 8\cdot 10^{11}$.

Furthermore, no trend is apparent looking at the last digit of the zeroes, being all digits (1, 3, 5, 7 and 9 only) equiprobable. Thus, we have shown that the distribution of zeroes is, if not "random"\footnote{How random are the zeroes in the sequence? Naively one would assume them to be as random as the prime number distribution.
	It is beyond the scope of this work to answer this question, but it would be worthy to study the randomness of this sequence (or any other $F_y(n)$ created with this scheme)
	by using advanced randomness tests (\citeauthor{links,marsaglia,tsang}).}, at least rather capricious.
On the other hand, the gaps between zeroes present a remarkably well defined trend, displaying an exponential decay according to the size, the length of the gap and the frequency with which it appears.
This important result echoes previous studies regarding gaps between primes, such as those presented in (\citeauthor{wolf2} (\citeyear{wolf2}), \citeauthor{szpiro1} (\citeyear{szpiro1}), \citeauthor{szpiro2} (\citeyear{szpiro2}), \citeauthor{ares, wolf4} (\citeyear{wolf4}), \citeauthor{ball}, \citeauthor{wolf5} (\citeyear{wolf5}), \citeauthor{garcia}).
Whether or not this will be true for any given interval is another open question.

Our theoretical predictions are largely based on intuition, but we provide qualitative support from our preliminary results. Although the behaviors reported here have been validated only for $n$ as large as $8\cdot 10^{11}$, it is reasonable to expect that they would be observed for any arbitrary $n$. However, mathematical proofs for our conjectures are likely to be extremely difficult since they are inextricably intertwined with the still mysterious distribution of prime numbers.

To conclude, the algorithm\footnote{The program is available under request.} of Jacob’s Ladder exhibits a remarkable feature despite its simplicity: Most of the zeroes are close to prime numbers since the smaller the gaps are, the more frequently they appear. This important result represents an unexpected correlation between the apparently chaotic sequence, the zeroes in $J(n)$, and the prime numbers distribution.

\appendixtitles{yes}

\appendix
\section{Additional conjectures}
\label{app1}

\vspace{6pt}

{\bf Conjecture II:}

The slope $\epsilon(n)$, of the Ladder is zero in the limit when $n$ goes to infinity.
\vspace{6pt}

\underline{Discussion}:

The Ladder could have a finite (although we conjecture this not to be case) amount of zeroes, so that after a certain number $X$, all points would be either above or below the $x$ axis, and hence the slope would be positive or negative in consequence.
Note that even if Conjecture I is false, the slope could decrease continuously when $n$ increases. And obviously if Conjecture I is true, then Conjecture II is more likely, but still not necessarily true.
\begin{wrapfigure}{R}{0.45\textwidth}
	\begin{center}
		\includegraphics[width=.9\linewidth]{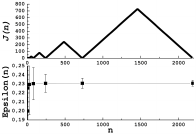}
		\caption{\emph{Upper panel:} A counterexample of a ladder presenting infinitely many zeroes, but where $\epsilon \neq 0$. \emph{Lower panel:} $\epsilon$ vs the number of points in the ladder.}
		\label{fig:counterattack}
	\end{center}
\end{wrapfigure}\leavevmode

It is important to note that an arbitrarily large amount of zeroes does not imply $\epsilon = 0$. This can be proved with a simple counterexample, as shown in Fig. \ref{fig:counterattack} (upper panel), where a simple ladder is constructed by stacking triangles, each one three times higher than the previous one. This ladder will have an infinite number of zeroes by construction, but as can be easily seen, the slope $\epsilon$ is not zero (The lower panel in Fig. \ref{fig:counterattack} shows the values of $\epsilon$ vs the number of points).

A clearer picture of the late behavior of Jacob's ladder can be ascertained from an analysis of the gaps and the frequency with which they appear, as we discuss in Sec. \ref{sec:Discussion}. Under closer inspection, a pattern with signs of predictability arises, allowing us to get an intuition on how the ladder is likely to extend towards infinity, much in the same way that a physicist would extrapolate statistically. To which extent one should trust that intuition is not a matter of debate for the present paper.

For these reasons, we argue that Conjecture II is true even if Conjecture I is not. Now a different, more complex question is how fast $\epsilon(n) \to 0$. The value of $\epsilon(n)$ is obtained by fitting the Ladder to $y=b\,x$ following a least squares approach
\be
b=\frac{n\sum x_i y_i -\sum x_i \sum y_i}{n\sum x_i^2-\left(\sum x_i\right)^2}\,.
\ee
Whether the slope $b$ tends to zero or not depends on the difference between the terms of the numerator, which is where the prime numbers define the Ladder through the $y_i$, which are unknown a priori (the $x_i$ values are known for all $i$). Many important papers have been published in the last decades - see (\citeauthor{wolf} (\citeyear{wolf}), \citeauthor{gamba,cramer,erdos,maynard1} (\citeyear{maynard1}), \citeauthor{zhang,guy,hardy2} (\citeyear{hardy2}), \citeauthor{polymath1,maynard3} (\citeyear{maynard3}), \citeauthor{stein,wolf2} (\citeyear{wolf2}), \citeauthor{wolf3} (\citeyear{wolf3}), \citeauthor{odlyzko,szpiro1} (\citeyear{szpiro1}), \citeauthor{szpiro2} (\citeyear{szpiro2}), \citeauthor{ares,wolf4} (\citeyear{wolf4}), \citeauthor{ball,wolf5} (\citeyear{wolf5}), \citeauthor{garcia}) to name a few - some of which could help or show the way to calculate certain upper (and/or lower) limits for $\epsilon(n)$ in the limit when $n$ tends to infinity. It seems to be a problem that could be attacked using methods of probabilistic number theory (\citeauthor{kubilius,kowalski}), however, it is beyond the scope of this paper to prove this conjecture.
\begin{figure}[!b]
	\centering
	\includegraphics[width=0.72\textwidth]{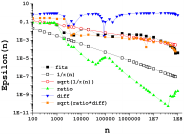} \\
	\caption{(Color online) Slope $\epsilon(n)$ of Jacob’s Ladder from 100 to 100,000,000 in logarithmic scale. Black full squares are the numerical values (obtained by fitting the Ladder to $y=b\,x$), while red empty circles, black empty squares and blue/green triangles are simple models. Orange squares are the values obtained averaging both green and blue points (see text).}
	\label{fig:sixth}
\end{figure}
 
A first naive idea could be to assume that the slope is of the same order of magnitude as $\frac{n/\pi(n)}{n}$, that is to say, $1/\pi(n)$, where $\pi(n)$ is the prime-counting function as usual (See Fig. \ref{fig:sixth}). It is to be mentioned that Fig. \ref{fig:sixth} plots the absolute value of b , but this value can be positive or negative. It goes without saying that the slope will depend not only on the number of primes, but also on their order, which is unknown a priori.\footnote{The slope is plotted up to $10^8$ due to computational limitations. A way of calculating the slope of billions of points would be to "sift" or "filter" them, that is to say, to take only 1 every 100 or every 1000, etc. The error, however, is rather difficult to determine.}.
The $1/\pi(n)$ model seems however to be an underestimate of $\epsilon(n)$, suggesting that $\sqrt{1/\pi(n)}$ could be a better approximation. Obviously one can always obtain a better fit by adding one more parameter to some function, in this case, a power to $1/\pi(n)$.

In any case, this simple expression fits reasonably well the $\epsilon(n)$ values for $n > 10,000$ (See Fig. \ref{fig:sixth}, red open circles).
It is also instructive to consider two very simple models, which are presented in blue (green) triangles. These are two clearly crude estimates, obtained by dividing the difference (ratio) between points up and down over $n$. We conjecture that the ratio $p$ will tend to 1 (see next section), so the slope calculated this way will tend to 0. It can be seen that the models badly over- and underestimate the value of $\epsilon(n)$ almost in the whole range, but interestingly, the geometric average (i.e. the square root of the product) of both functions seems to be a valid ballpark estimate.

\vspace{12pt}

{\bf Conjecture III-A}

If Conjecture I is true, then the area below the upper part of the Ladder is equal to the area above the lower part of the Ladder when $n\to\infty$, or more precisely, the ratio between both areas tends to 1 in that limit.

\be
\lim^{}_{n\to\infty} \frac{A_{\text{pos}}(n)}{A_{\text{neg}}(n)}=1
\ee

\underline{Discussion}:

It is natural to think that if no particular order is found, then the ratio of both areas in the limit will be 1. This is similar, but not equivalent to say that the number of points above and below $y=0$ (positive and negative respectively, or just $C_{\text{pos}}$ and $C_{\text{neg}}$) will be the same for $n\to\infty$, so that their ratio will tend to 1. Thus, a similar conjecture can be formulated:

\vspace{12pt}

{\bf Conjecture III-B}

If Conjecture I is true, the ratio between numbers of points above and below $y=0$ will tend to 1.

\be
\lim^{}_{n\to\infty} \frac{C_{\text{pos}}(n)}{C_{\text{neg}}(n)}=1
\ee

\underline{Discussion}:
\begin{figure}[!t]
	\centering
	\includegraphics[width=0.72\textwidth]{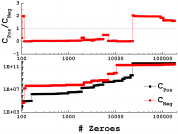} \\
	\caption{(Color online) \emph{Bottom:} Number of positive (black squares) and negative (red circles) points in Jacob’s Ladder up to  $8\cdot 10^{11}$. \emph{Top:} Ratio between the number of positive and negative points up to $8\cdot 10^{11}$. The blue dashed line marks $y=1$, the conjectured asymptotic ratio. Note the logarithmic scale in both axes.}
	\label{fig:seventh}
\end{figure}

The area depends not only on the number of points up or down but also on the ordinates of those points, so it is straightforward to see that Conjecture III-A can be true and Conjecture III-B false, or the other way around. However, we presume that in the limit $n\to\infty$ both the area and the ratio $C_{\text{pos}}/C_{\text{neg}}$ will tend to 1. Fig. \ref{fig:seventh} (bottom) shows the number of positive (black) and negative (red) points in Jacob’s Ladder up to $8\cdot 10^{11}$ (plotted against the number of zeroes to show it as a continuum). It is to be mentioned as well that there is a strong possibility that these limits do not exist.

Importantly, note that Conjecture III (A or B) is not equivalent to Conjecture I. If Conjecture III (A or B) is true then Conjecture I is true but it could be that Conjecture I is true and Conjecture III (A or B) is not\footnote{Conjecture III (A and B) could be stated as Conjecture I, and hence we could talk of Conjecture I (infinitely many zeroes) as a corollary, given that Conjecture III (A, B or both) demonstrates that the number of zeroes tends to infinity, being the result a corollary of Conjecture I. However we consider more natural, at a first glance of Jacob’s ladder, to think of the number of possible cuts (zeroes) when $n\to\infty$, so we decided to keep them in this order.}.

Fig. \ref{fig:seventh} (top) shows the ratio between the number of positive and negative points in Jacob’s Ladder (presented in Fig. \ref{fig:sixth}) up to 100 million. As can be seen, some approximate periodicity over a decreasing behaviour seems apparent. Nonetheless, we believe this to be a misleading result so that the ratio, whenever $n$ is sufficiently large, will tend to 1, and the approximate periodicity in Fig. \ref{fig:seventh} (top) will disappear.

\vspace{1cm}

\authorcontributions{Conceptualization, A.F. and R.M.; methodology, R.M. and D.F.; software, R.M. and D.F.; validation, A.F., R.M. and D.F.; formal analysis, A.F.; investigation, A.F.; writing--original draft preparation, A.F.; writing--review and editing, D.F.}

\funding{This research was supported by the Ministry of Education, Youth and Sport of the Czech Republic under project CZ.02.1.01/0.0/0.0/15\_003/0000464, and by the Icelandic Research Fund, under contract 163419-051. AF acknowledges the use of the Anselm and Salomon HPC clusters administered by IT4Innovations National Supercomputing Center, University of Ostrava. DF is grateful for the hospitality of the Technical University of Vienna.}

\conflictsofinterest{On behalf of all authors, the corresponding author states that there is no conflict of interest.} 

\reftitle{References}


\begin{thebibliography}{999}

\bibitem[Ingham(1932)]{ingham}
Ingham, A.E. The Distribution of Prime Numbers. In {\em Cambridge Mathematical Tracts} {\bf 1932}, {\em 30}; Cambridge University Press.

\bibitem[Littlewood(1914)]{littlewood}
Littlewood, J.E. Sur les distribution des nombres premiers. {\em Comptes Rendus Acad. Sci. Paris} {bf 1914} {\em 158}, 1869--1872; Littlewood, J.E. Sur la distribution des nombres premiers. {\em Comptes Rendus} {\bf 1914}, {\em v. 158}, pp. 1869-1872.

\bibitem[Anglin(1995)]{anglin}
Anglin, W.S. {\em The Queen of Mathematics: An Introduction to Number Theory}; Dordrecht, Netherlands: Kluwer, 1995.

\bibitem[Apostol(1976)]{apostol}
Apostol, T.M. {\em Introduction to Analytic Number Theory}; New York: Springer-Verlag, 1976.

\bibitem[Hardy(1979)]{hardy}
Hardy, G.H. and Wright, W.M. Unsolved Problems Concerning Primes. In {\em An Introduction to the Theory of Numbers}; 5th Ed.; Oxford University Press: Oxford, UK, 1979; pp. 19 and 415--416.

\bibitem[Shanks(1985)]{shanks}
Shanks, D. {\em Solved and Unsolved Problems in Number Theory}; New York: Chelsea, 1985; pp. 30-31 and 222.

\bibitem[Stallings(1999)]{stallings}
Stallings, M. {\em Cryptography and Network Security: Principles and Practice}; Prentice-Hall; New Jersey, 1999.

\bibitem[Goles(2001)]{goles}
Goles, E., Schulz, O. and Markus, M. Prime number selection of cycles in a predator-prey model. {\em Complexity} {\bf 2001}, {\em 6}, 33–38.

\bibitem[Toha(1999)]{toha}
Toha, J. and Soto, M.A. Biochemical identification of prime numbers. {\em Medical Hypothesis} {bf 1999}, {\em 53}, 361.

\bibitem[Berry(1993)]{berry}
Berry, M.V. Quantum chaology, prime numbers and Riemann’s zeta function. {\em Inst. Phys. Conf. Ser.} {\bf 1993}, {\em 133}, 133–134.

\bibitem[Sakhr(2003)]{sakhr}
Sakhr, J., Bhaduri, R.K. and Van Zyl, B.P. Zeta function zeros, powers of primes, and quantum chaos. {\em Phys. Rev. E} {\bf 2003}, {\em 68}, 026206.

\bibitem[Wolf(1989)]{wolf}
Wolf, M. Multifractality of prime numbers. {\em Physica A} {\bf 1989}, {\em 160}, 24.

\bibitem[Gamba(1990)]{gamba}
Gamba, Z., Hernando J. and Romanelli, L. Are prime numbers regularly ordered?. {\em Phys. Lett. A} {\bf 1990}, {\em 145}, 106.

\bibitem[Saouter(1998)]{saouter}
Saouter, Y. Checking the odd Goldbach conjecture up to $10^{20}$. {\em Math. Comput.} {\bf 1998}, {\em 67}, 863-866.

\bibitem[Richstein(2001)]{richstein}
Richstein, J. Verifying the Goldbach conjecture up to $4\cdot 10^{14}$. {\em Math. Comput.} {\bf 2001}, {\em 70}, 1745-1750.

\bibitem[Granville(1989)]{granville}
Granville, A., Van de Lune, J. and Te Riele, H.J.J. Checking the Goldbach conjecture on a vector computer. {\em Proceedings of NATO Advanced Study Institute} {\bf 1989}, {\em 1988}, 423–433.

\bibitem[Helfgott(2013)]{helfgott}
Helfgott, H.A. The ternary Goldbach conjecture is true. \href{https://arxiv.org/abs/1312.7748}{ArXiv[1312.7748]}

\bibitem[Borwein(2003)]{borwein}
Borwein, J. and Bailey, D. {\em Mathematics by Experiment: Plausible Reasoning in the $21^{\text{st}}$ Century}; A. K. Peters Co.; Wellesley, MA, 2003; p. 64.

\bibitem[Cramer(1936)]{cramer}
Cram\'er, H. On the order of magnitude of the difference between consecutive prime numbers. {\em Acta Arithmetica} {\bf 1936}, {\em 2}, 23–46.

\bibitem[Erdos(1935)]{erdos}
Erd\'os, P. On the difference of consecutive primes. {\em Q. J. Math. Oxford Ser.} {\bf 1935}, {\em 6},124–128.

\bibitem[Maynard(2016)]{maynard1}
Maynard, J. Large gaps between primes. {\em Ann. of Math.} {\bf 2016}, {\em 183}, iss. 3, pp. 915-933; Maynard, J. Dense clusters of primes in subsets. {\em Compositio Math.} {\bf 2016}, {\em 152}, no. 7, pp. 1517-1554.

\bibitem[Zhang(2014)]{zhang}
Zhang, Y. Bounded gaps between primes. {\em Ann. Of Math Second Ser.} {\bf 2014}, {\em 179(3)}, 1121–1174.

\bibitem[Guy(1994)]{guy}
Guy, R. {\em Unsolved Problems in Number Theory}; 2nd Ed., Springer, 1994; p. vii.

\bibitem[Hardy(1923)]{hardy2}
Hardy, G.H. and Littlewood, J.E. Some problems of 'Partitio numerorum' III: On the expression of a number as a sum of primes. {\em Acta Math.} {\bf 1923}, {\em 44}, 1-70.

\bibitem[Polymath(2014)]{polymath1}
Polymath, D.H.J. New equidistribution estimates of Zhang type, and bounded gaps between primes. {\em Algebra Number Theory} {\bf 2014}, {\em 8}, 2067-2199; Polymath, D.H.J. Variants of the Selberg sieve, and bounded intervals containing many primes. {\em Research in the Mathematical Sciences} {\bf 2014}, {\em 1:12}.

\bibitem[Maynard(2015)]{maynard3}
Maynard, J. Small gaps between primes. {\em Ann. of Math.} {\bf 2015}, {\em 181}, iss. 1, pp. 383-413.

\bibitem[Stein(1964)]{stein}
Stein, M.L., Ulam, S.M. and Wells, M.B. A visual display of some properties of the distribution of primes. {\em The American Mathematical Monthly} {\bf 1964}, {\em 71}, No. 5, pp. 516-520.

\bibitem[Orlowski(2016)]{orlowski}
Chmielewski, L.J. and Orlowski, A. Finding Line Segments in the Ulam Square with the Hough Transform. {\em Computer Vision and Graphics - Lecture Notes in Computer Science} {\bf 2016}, {\em 9972}.

\bibitem[Hull(1962)]{hull}
Hull, T.E. and Dobell, A.R. Random Number Generators. {\em SIAM Review} {\bf 1962}, {\em 4}, 230–254.

\bibitem[Billingsley(1973)]{billingsley}
Billingsley, P. Prime numbers and Brownian motion. {\em Amer. Math. Monthly} {\bf 1973}, {\em 80}, 1099-1115.

\bibitem[Peng(1992)]{peng}
Peng, C.-K., Buldyrev, S.V., Goldberger, A.L., Havlin, S., Sciortino, F., Simons, M. and Stanley, H.E. Long-range correlations in nucleotide sequences. {\em Nature} {\bf 1992}, {\em 356}, 168.

\bibitem[Wolf(1998)]{wolf2}
Wolf, M. Random walk on the prime numbers. {\em Physica A} {\bf 1998}, {\em 250}, 335-344.

\bibitem[Wolf(1996)]{wolf3}
Wolf, M. Unexpected regularities in the distribution of prime numbers. {\em Proceedings of the $8^{\text{th}}$ Joint EPS-APS International Conference, Krakow} {\bf 1996}, pp. 361–367.

\bibitem[Odlyzko(1999)]{odlyzko}
Odlyzko, A., Rubinsten, M. and Wolf, M. Jumping champions. {\em Exp. Math.} {\bf 1999}, {\em 8}, 107–118.

\bibitem[Szpiro(2004)]{szpiro1}
Szpiro, G.G. The gaps between the gaps: some patterns in the prime number sequence. {\em Physica A} {\bf 2004}, {\em 341}, 607–617.

\bibitem[Szpiro(2007)]{szpiro2}
Szpiro, G.G. Peaks and gaps: Spectral analysis of the intervals between prime numbers. {\em Physica A} {\bf 2007}, {\em 384}, 291–296.

\bibitem[Ares(2006)]{ares}
Ares, S. and Castro, M. Hidden structure in the randomness of the prime number sequence? {\em Physica A}, {\bf 2006}, {\em 360}, 285–296.

\bibitem[Wolf(1997)]{wolf4}
Wolf, M. 1/f noise in the distribution of primes. {\em Physica A} {\bf 1997}, {\em 241}, 439–499.

\bibitem[Ball(2003)]{ball}
Ball, P. Prime numbers not so random? {\em Nature Science Update} {\bf 2003}.

\bibitem[Wolf(2014)]{wolf5}
Wolf, M. Nearest-neighbor-spacing distribution of prime numbers and quantum chaos. {\em Phys Rev E} {\bf 2014}, {\em 89}, 022922.

\bibitem[Garcia-Perez(2014)]{garcia}
Garc\'ia-Perez, G., Serrano, M.A. and Bogu\~na, M. Complex architecture of primes and natural numbers. {\em Phys Rev E} {\bf 2014}, {\em 90}, 022806.

\bibitem[Benford(1938)]{benford}
Benford, F. The law of anomalous numbers. {\em Proc. Am. Philos. Soc.} {\bf 1938}, {\em 78}, 551–572.

\bibitem[Klafter(2000)]{klafter}
Metzler, R. and Klafter, J. The random walk's guide to anomalous diffusion: A fractional dynamics approach. {\em Phys. Rep.} {\bf 2000}, {\em 339}, 1-77.

\bibitem[Newcomb(1881)]{newcomb}
Newcomb, S. Note on the frequency of use of the different digits in natural numbers, {\em Am. J. Math.} {\bf 1881}, {\em 4}, 39–40.

\bibitem[Draconis(1977)]{draconis}
Diaconis, P. The distribution of leading digits and uniform distribution mod 1. {\em Ann. Probab.} {\bf 1977}, {\em 5}, 72–81.

\bibitem[Raimi(1976)]{raimi}
Raimi, R.A. The first digit problem. {\em Am. Math. Mon.} {\bf 1976}, {\em 83}, 521–538.

\bibitem[Takloo-Bighash(2006)]{takloo}
Miller, S.J. and Takloo-Bighash, R. {\em An invitation to modern number theory}; NJ: Princeton University Press; Princeton, 2006.

\bibitem[Dickson(1919)]{new1}
Dickson, L.E. {\em History of the Theory of Numbers}; Washington, 1919.

\bibitem[Tschebyschef(1853)]{new2}
Tschebyschef, P.L. Lettre théorème relatif aux nombres premiers de M. le professeur Tche'bychev à M. Fuss sur un nouveaux contenus dans les formes $4n \pm 1$ et $4n \pm 3$. {\em Bull, de la Classe Phys. de l'Acad. Imp. des Sciences, St. Petersburg} {\bf 1853}, {\em v. 11}, p. 208.

\bibitem[Phragmen(1892)]{new3}
Phragmen, E. Sur le logarithme integral et la fonction $f(x)$ de Riemann. {\em Ofversigt af Kongl. Vetenskaps-Akademiens Forhandlingar} {\bf 1891-92}, {\em 48}, 721 —744.

\bibitem[Landau(1953)]{new4}
Landau, E. {\em Handbuch der Lehre von der Verteilung der Primzahlen}; New York, 1953.

\bibitem[Luque(2009)]{bartolo}
Luque, B. and Lacasa, L. The first-digit frequencies of prime numbers and Riemann zeta zeros. {\em Proc. R. Soc. A} {\bf 2009}, {\em 465}, 2197–2216.

\bibitem[OEIS(2019)]{link}
Entry in the On-Line Encyclopedia of Integer Sequences: \href{https://oeis.org/A065358}{https://oeis.org/A065358}

\bibitem[Cacert(2019)]{links}
Cacert research lab's: \href{http://www.cacert.at/random/}{http://www.cacert.at/random/} and the {\bf ent} program: \href{http://www.fourmilab.ch/random/}{http://www.fourmilab.ch/random/} 

\bibitem[Marsaglia(1985)]{marsaglia}
Marsaglia, G. {\em The Marsaglia Random Number CDROM, with The Diehard Battery of Tests of Randomness} {\bf 1985}, 
produced at Florida State University under a grant from The National Science Foundation, 
available at \href{http://www.cs.hku.hk/diehard/}{http://www.cs.hku.hk/diehard/} or an earlier version at \href{http://www.stat.fsu.edu/pub/diehard}{http://www.stat.fsu.edu/pub/diehard}.

\bibitem[Tsang(2002)]{tsang}
Marsaglia, G. and Tsang, W.W. Some difficult-to-pass tests of randomness. {\em Journal of Statistical Software} {\bf 2002}, {\em 7}, issue 3.

\bibitem[Kubilius(1962)]{kubilius}
Kubilius, J. Probabilistic methods in the theory of numbers. {\em Translations of mathematical monographs} {\bf 1962}, {\em 11}; RI: American Mathematical Society, Providence, 1962.

\bibitem[Kowalski(2019)]{kowalski}
Kowalski, E. Arithmetic Randonnée: an introduction to probabilistic number theory. {\em ETH Zürich Lecture Notes} {\bf 2019}, \href{www.math.ethz.ch/~kowalski/probabilistic-number-theory.pdf}{www.math.ethz.ch/~kowalski/probabilistic-number-theory.pdf}.

\end{thebibliography}
\end{document}